\documentclass{article}
\usepackage[utf8]{inputenc}
\usepackage{amsmath}
\usepackage{amssymb}
\usepackage{amsfonts}
\usepackage[russian]{babel}

\title{Schottky uniformization and Bloch-Wigner dilogarithm on higher genus curves}
\author{Ilyas Bayramov}
\date{November 2020}

\begin{document}

\maketitle

\section{Motivation}

The integral that defines the dilogarithm first appears in one of the letters from Leibniz to Johan Bernoulli in 1696 (Maximon)[7], and is equal to $$Li_2(z)=-\int_{0}^z \frac{\text{ln}(1-u)}{u}du.$$ This can be rewritten as an iterated integral $$Li_2(z)=\int_0^z \frac{dt}{t}\int_0^t\frac{du}{1-u}.$$ Since its introduction, it has been studied by innumerable mathematicians, among which we can name Euler who defined it as the power series $$Li_2(z)=\sum^\infty_{k=1}\frac{z^k}{k^2},$$ Abel, Lobachevsky, Kummer, and Ramanujan, among others, and appears in a variety of fields both in pure mathematics-such as number theory, algebraic geometry, hyperbolic geometry and combinatorics-and applied, in problems related to electrical networks, radiation, statistical mechanics of ideal gases, and quantum electrodynamics in any calculation of higher-order processes such as vacuum polarization and radiative corrections. It is the ur-example of an iterated integral.

The single-valued version of dilogarithm was defined by Bloch and Wigner [2] as  $\text{D}_2 (z) = \text{Im (Li}_2 (z) )+\text{arg}(1-z)\text{log}|z|.$ The result was generalized by Ramakrishnan to higher polylogarithms which are defined iteratively as $\text{Li}_n(z)=\int_0^z\frac{\text{Li}_{n-1}(t)}{t}dt$, as follows: let $\text{L}_m(z)=\sum^m_{j=1}\frac{(-\text{log}|z|)^{m-j}}{(m-j)!}\text{Li}_j(z),$ then $\text{D}_m=\text{Im}(L_m),$ if $m$ is even, and $\text{Re}(L_m)+\frac{(\text{log}(|x|))^m}{2m!}$ if $m$ is odd. Bloch has defined the elliptic version of $\text{D}_2$, $\text{D}_2(q,x)=\sum_{k\in\mathbb{Z}} \text{D}_2(k,x), |x|>0, 0<|q|<1,$ and proved its convergence, Zagier has generalized it to higher polylogarithms (Zagier [11]), and Levin has generalized both of these results to the elliptic version of $\text{Li}_n$ for all $n$ (Levin [5]).

In this paper, we generalize Bloch's result to higher genus curves.\\

\textbf{Theorem 1:} Let $\text{D}:=\text{D}_2$ be as above. Then for any Riemann surface $S$ there exists a Schottky uniformization by a group $\Gamma_S$ such that the Poincare series of $\text{D}$, $\text{D}_{\Gamma}(z)=\sum_{\gamma \in \Gamma}\gamma'(z)\text{D}(\gamma z)$ converges, and is a well-defined single-valued function on $S.$\\

This result exploits Patterson-Sullivan theory and recent advancements in Schottky uniformization that are based on it.

\section{Acknowledgements}

I thank Andrey Levin in his capacity as my academic advisor. I also thank Alexey Glutsyuk for many valuable discussions. Finally, I would like to thank Renat Gontsov for refereeing this paper at its Masters thesis stage.

\section{Schottky uniformization}

\textbf{Koebe's Retrosection theorem.}[6] Any Riemann surface can be parameterized by a Schottky group.\\
\\
\textbf{Proof.} Consider the Riemann surface $S$ of genus $g.$ We can chose a collection $\{Q_i\}_{1,...,g}$ of simple closed curves on $S$, such that the homology classes they represent are linearly independent and do not intersect (that is, the canonical symplectic form on $H_1(S,\mathbb{Z})$ is zero on the pairs of their homology representatives.) Cutting along those curves, we obtain a topological sphere with the interior of $2g$ Jordan curves removed. By uniformization theorem, this sphere is in fact biholomorphically equivalent to the Riemann sphere with the interiors of $2g$ Jordan curves removed. Let's call it $A.$ Then the gluing maps are in fact Mobius transformations, and in fact hyperbolic, so the group generated by them is Schottky.\\

Now, a question one might ask is how large is the lower bound on the Hausdorff dimension of the limit set of the Schottky groups that uniformize a given Riemann surface. This has been a long-standing problem, finally settled by Hou in 2016.\\

\textbf{Theorem(Hou)}[4]: Any Riemann surface can be uniformized by a Schottky group with Hausdorff dimension less than $1$.\\

 The idea of Hou's proof is to construct a measure on $W(\Gamma)\times \Lambda(\Gamma),$ where $W(\Gamma)$ is the collection of generating sets of $\Gamma$ that will detect the action of the mapping class group (the Patterson-Sullivan measure is well-suited for it, because of its quasi-invariance property).\\

\textbf{Definition:} The \textbf{Schottky space} is the space that parameterizes all Schottky groups up to conjugation in $PSL(2,\mathbb{C}).$\\

We can easily see that this space is parameterized by the coordinates of the fixed points of the hyperbolic transformations (there are $2g$ such extended complex numbers) and the highest eigenvalue of these transformations ($g$ such complex numbers,) up to conjugation (the Mobius group is $3$-transitive, so we can fix $3$ of the first $2g$ numbers to be $0,1,\infty$), so the overall dimension is $3g-3.$\\

Now, from Koebe's uniformization, a Schottky group that uniformizes a given Riemann surface is defined by mapping $g$ of the generators of a basis of $H_1(S,\mathbb{Z})$ (the ones that do not intersect) to a free subgroup of $SL(2,\mathbb{C}).$\\

We shall consider the space of Riemann surfaces up to diffeomorphisms that are isotopic to the identity. This is the \textbf{Teichmuller space}. The group $\pi_0(Diff(S))$ is called the mapping class group, and it acts discontinuously on the Teichmuller space. The quotient of the Teicmuller space under this action is precisely the \textbf{moduli space}. We can further choose the homology marking, that is, fix a symplectic basis of the homology, on the Riemann surface. In this way, we obtain the \textbf{Torelli space.} The elements of the mapping class group that fix the homology classes pointwise are precisely the elements of the \textbf{Torelli group.}\\

Now, the Torelli space is a covering of the Schottky space, by forgetting half of the homological marking, but also up to the outer automorphisms of the Schottky group. The Schottky space covers the moduli space of Riemann surfaces, that is, the space of complex structures on a topological surfaces up to diffeomorphisms.\\

Precisely, the outer automorphisms of the Schottky group are generated Nielsen transformations-the non-commutative analogues of the usual elementary row operations, that replace a given generator by $x^{-1}$, permute two generators, replace a generator by its product with another one, and permute the generators cyclically-that are restricted only to the complement of the cut system. Then, the upper bound is established on the Hausdorff dimension of the limit set by taking the average of distortions over such transformations.\\

Using quasi-invariance of Patterson-Sullivan measure Nayatani managed to construct a Poincare metric on the domain of discontinuity $\Omega(\Gamma)$ of the action of a Kleinian group $\Gamma$ on $\mathbb{S}^2$:\\

\textbf{Theorem (Nayatani)}[8]: Let $\phi: \mathbb{S}^2\times\mathbb{S}^2\to \mathbb{R}, \phi(x,y)=1/2|x-y|^2=1-\text{cos} r,$ $r-$ the geodesic distance with respect to the induced metric $g_0$, $\mu-$ a Patterson-Sullivan measure with unit total mass. Then $$g_{\mu}=(\int_{\Lambda(\Gamma)}\phi^{-\delta}(\cdot,y)d\mu(y))^{2/\delta}g_0$$ is an invariant conformally flat metric on $\Omega(\Gamma),$ where $\delta<1$.\\
\textbf{Proof:} $$\int_{\Lambda(\Gamma)}\phi^{-\delta}(\gamma x,y)d\mu(y)=\int_{\Lambda(\Gamma)}\phi^{-\delta}(\gamma x,\gamma y)d\mu(\gamma y)=|\gamma'(x)|^{-\delta}\int_{\Lambda(\Gamma)}\phi^{-\delta}(x,y)d\mu(y),$$ since $\phi(\gamma x, \gamma y)=|\gamma'(x)||\gamma'(y)|\phi(x,y),$ and $\gamma^{*}\mu(y)=|\gamma'(y)|^{\delta}\mu(y)$ so the metric $g_{\mu}$ is clearly invariant, and the conformal flatness follows from the definition of the metric as a multiple of $g_0,$ which is also conformally flat. We can take $\delta<1$ due to Hou.\\

The existence of such a metric and its asymptotic growth near $\Lambda(\Gamma)$ allows us to explicitly verify the following criterion, derived by Bers, on the convergence of the Poincare series of $\Gamma$ of a continuous function on a domain of discontinuity of $\Gamma$.\\

\textbf{Remark (Bers)}[1]: Let $\Delta$ be a domain in $\mathbb{C},$ let $\Gamma$ act on $\Delta$ properly discontinuously. Then if $\rho$ is a positive continuous function on $\Delta$ such that $\rho(\gamma(\zeta))\text{jac}(\gamma)=\rho,$ for $\gamma\in \Gamma$ $\Phi$ a continuous function on $\Delta,$ the Poincare series of $\Phi$ converges iff $\int_{\Delta} \rho(\zeta)^2\Phi(\zeta)d\zeta$ does.\\

\textbf{Proof:} This is the usual integral test.\\

\textbf{Proof (Theorem 1):}

In our case, $\rho$ is exactly the function $(\int_{\Lambda(\Gamma)}\phi^{-\delta}(\gamma \cdot,y)d\mu(y))^{1/\delta}$ in spherical coordinates, and $\int_{\Lambda(\Gamma)}\phi^{-\delta}(\gamma \cdot,y)d\mu(y)\sim  x^{\delta(\delta-1)},$ where $x=|\zeta-y|,$ and $y$ is finite, and $|\zeta|^{-\delta^2}$, as $\zeta\to\infty$ because of stereographic transform; using the asymptotics of Bloch-Wigner function [11] ($D(|x|)\sim  |x|\text{log}|x|$ as $x\to 0$) and its modular property $D(1/x)=(-1)D(x)$, we see that $(\rho(x))^{2/\delta}D(x)\sim  x^{2\delta-1}\text{log}(x),$ so the integral in question converges, as $\delta<1;$ on other hand, as $D(x)\sim D(1/x)),$ we have that lim$_{x\to\infty}\rho^{2/\delta}(x)D(x)=$lim$_{x\to \infty} |x|^{-1}|x|^{-2\delta}\text{log}|x|,$ and the power of $x$ in the denominator is greater than $1,$ since the Hausdorff dimension of the limit set cannot be $0$ unless the Schottky group is cyclic, so the integral from Bers's remark converges. Hence, the Poincare series of the Bloch-Wigner function is convergent.

\newpage
\section{References}

\begin{description}\item[]
\leavevmode \\
1. Lipman Bers, Automorphic Forms for Schottky Groups. In "Surveys in Applied Mathematics: Essays Dedicated to S.M. Ulam", ed. by N. Metropolis, S. Orszag and G.-C. Rota. Academic Press, 1976.
\end{description}

\begin{description}\item[]
\leavevmode \\
2. Spencer Bloch, "Higher Regulators, Algebraic K-theory and Zeta Functions of Elliptic Curves". AMS, 2000.\\
\end{description}

\begin{description}\item[]
\leavevmode \\
3. Pierre De La Harpe, "Topics in Geometric Group Theory". University of Chicago Press, 2000.\\
\end{description}

\begin{description}\item[]
\leavevmode \\
4. Yong Hou, "On smooth moduli space of Riemann surfaces". \\https://arxiv.org/abs/1610.03132.\\
\end{description}

\begin{description}\item[]
\leavevmode \\
5. Andrey Mikhailovich Levin, "Elliptic polylogarithms: An analytic theory". In Compositio Mathematica, CUP, 1997.\\
\end{description}

\begin{description}\item[]
\leavevmode \\
6.Bernard Maskit, "Kleinian Groups". Springer, 2004.\\

\end{description}

\begin{description}\item[]
\leavevmode \\
7. Leonard Maximon, "The Dilogarithm Function for Complex Argument". In "Proceedings: Mathematical, Physical and Engineering Sciences". The Royal Society, 2003.\\

\end{description}

\begin{description}\item[]
\leavevmode \\
8. Shin Nayatani, "Patterson-Sullivan measure and conformally flat metrics". In "Mathematische Zeitschrift". Springer, 1997.\\

\end{description}

\begin{description}\item[]
\leavevmode \\
9. Peter Nicholls, "The Ergodic Theory of Discrete Groups".Cambridge University Press, 1989.\\

\end{description}

\begin{description}\item[]
\leavevmode \\
10. Jean-Francois Quint, "An Overview of Patterson-Sullivan Theory".\\ https://www.math.u-bordeaux.fr/~jquint/publications/courszurich.pdf.\\

\end{description}

\begin{description}\item[]
\leavevmode \\
11. Don Zagier, "The Bloch-Wigner-Ramakrishnan polylogarithm function". In "Mathematische Annalen", 1990.

\end{description}

\end{document}